# Derivations and Iterated Skew Polynomial Rings

Michael Gr. Voskoglou

**Abstract-** Two are the objectives of the present paper. First we study properties of a differentially simple commutative ring R with respect to a set D of derivations of R. Among the others we investigate the relation between the D-simplicity of R and that of the local ring $R_P$ with respect to a prime ideal P of R and we prove a criterion about the D- simplicity of R in case where R is a 1-dimensional (Krull dimension) finitely generated algebra over a field of characteristic zero and D is a singleton set. The above criterion was quoted without proof in an earlier paper of the author. Second we construct a special class of iterated skew polynomial rings defined with respect to finite sets of derivations of a ring R (not necessarily commutative) commuting to each other. The important thing in this class is that, if R is a commutative ring, then its differential simplicity is the necessary and sufficient condition for the simplicity of the corresponding skew polynomial ring.

## I. DIFFERENTIALLY SIMPLE COMMUTATIVE RINGS

**A**LL the rings considered in the first two sections of this paper are commutative with identity and all the fields are of characteristic zero, unless it is otherwise stated. A *local ring* is understood to be a Noetherian ring with a unique maximal ideal; if R is not Noetherian, then we call it a *quasi-local ring*. A Noetherian ring with finitely many maximal ideals is called a *semi-local ring*. For special facts on commutative algebra we refer freely to [2],[13], and [21], while for the concepts of Algebraic Geometry used in the paper we refer freely to [4].
Let R be a ring (not necessarily commutative) and let D be a set of derivations of R. Then an ideal I of R is called a *D-ideal*, if $d(I) \subseteq I$ for all derivations d in D , and R is called a *D-simple ring* if it has no non trivial D-ideals. When D is a singleton set, say D={d} , then, in order to simplify our notation, I is called a *d-ideal* and R is called a *d-simple* ring respectively. In general R is called a *differentially simple ring,* if there exists at least one set D of derivations of R, such that R is a D-simple ring.
Every D-simple ring R contains the field
$F= C(R) \cap [\cap_{d \in D} ker\ d]$, where $C(R)$ denotes the center of R, and therefore, if R is of characteristic zero, then it contains the field of rational numbers.
Non commutative differentially simple rings exist in abundance, e.g. every simple ring R is D-simple for any set D of derivations of R, and that is why our concern is turned to commutative rings only. Some characteristic examples of commutative D-simple rings, where D is a non singleton set of derivations, are given below.

**Example 1.1**: The polynomial ring R=k[$x_1, x_2, \ldots, x_n$] over a field k is a D- simple ring, where
D={ $\frac{\partial}{\partial x_1}, \frac{\partial}{\partial x_2}, \ldots \frac{\partial}{\partial x_n}$ }.
*Proof:* If I is a non zero D-ideal of R and f is in I, we can write f= $\sum_{i=0}^{k} f_i x_1^i$ , with $f_i$ in $k[x_2, \ldots, x_n]$ for each i. Then
$\frac{\partial^k f}{\partial x_1^k}$ =k!$f_k$ is a non zero element of
$I \cap k[x_2, \ldots, x_n]$. Repeating the same argument for $f_k$ and keep going in the same way one finds, after n at most steps, that $I \cap k \neq \{0\}$. Thus I=R.-

**Example 1.2:** The power series ring R=$k[[x_1, x_2, \ldots, x_n]]$ over a field k is a D-simple ring, where D is as in example 1.1 .
The proof is the same as above.-

**Example 1.3:** Let R= $\frac{IR[x, y, z]}{(x^2 + y^2 + z^2 - 1)}$ be the coordinate ring of the real sphere (IR denotes the field of the real numbers). Consider the IR-derivations $d_1$ and $d_2$ of the polynomial ring IR[x,y,z] defined by $d_1$ : $x \rightarrow y + z$, $y \rightarrow z - x,\ z \rightarrow -x - y,$ and $d_2$ : $x \rightarrow y + 2z$, $y \rightarrow xyz - x,\ z \rightarrow -xy^2 - 2x$. Then, since $d_i(x^2 + y^2 + z^2 - 1)$=0, for i=1,2 , $d_i$ induces an IR-derivation of R, denoted also by $d_i$.
Set D = {$d_1$ , $d_2$ }, then R is a D-simple ring (cf. [38], Lemma 3.1).-
The following example concerns a non trivial case of a differentially simple ring of prime characteristic:

**Example 1.4:** Let k be a field of prime characteristic, say p, and let R and D be as in example 1.1. Then I=($x_1^p, x_2^p, \ldots, x_n^p$) is obviously a D-ideal of R and therefore each $\frac{\partial}{\partial x_i}$, i=1,2,…,n , induces a derivation, say $d_i$, of the factor ring R/I by $d_i(f+I) = \dfrac{\partial f}{\partial x_i} + I$, for all f in R.

Set $D^* = \{d_1, d_2, \ldots, d_n\}$, then R/I is a $D^*$-simple ring.

*Proof:* If A is a non zero $D^*$-ideal of R/I, then $A' = \{f \in k[x_1, x_2, \ldots, x_n] : f + I \in A\}$ is obviously a D-ideal of $k[x_1, x_2, \ldots, x_n]$ containing properly I. It becomes evident that there exists f in A′ with all its terms of the form $c x_1^{m_1} x_2^{m_2} \ldots x_n^{m_n}$, with c in k and $0 \le m_i < p$, for each i=1, 2,…,n. Let m be the greatest integer which appears as exponent of $x_1$ in the terms of f. Without loss of the generality we may assume that $m \neq 0$. Then

$$0 = \frac{\partial^m f}{\partial x_1^m} = m! f_1(x_2, \ldots, x_n) \text{ is in } A'. \text{ If } f_1 \text{ is not in k, we}$$

repeat the same argument for $f_1$ and we keep going in the same way until we find, after a finite number of steps, that $A' \cap k \neq \{0\}$. Thus A = R/I and the result follows.-

There is no general criterion known to decide whether or not a given ring R is differentially simple. However it seems that there is a connection between the differential simplicity and the Krull dimension of R (denoted by dim R). The following result demonstrates this connection, when R is of prime characteristic:

**Theorem 1.5:** Let R be a ring of prime characteristic, say p, and let D be a set of derivations of R, such that R is D-simple. Then R is a 0-dimensional, quasi-local ring.

*Proof:* Let M be a maximal ideal of R and let I be the ideal of R generated by the set $\{m^p : m \in M\}$. Then, since R is of characteristic p, I is a proper D-ideal of R, therefore the D-simplicity of R implies that I=(0). Thus M is contained in the nil radical, say N, of R and therefore M=N. Let now P be a prime ideal of R contained in M. Then, since N is equal to the intersection of all prime ideals of R and M=N, we get that M=P. Thus N is the unique prime ideal of R and this proves the theorem.-

As an immediate consequence of the above theorem, if R is a domain, then R is a field (since M=N=(0)) and therefore the interest is turned mainly to rings of characteristic zero. In this case it is well known that a differentially simple ring is always a domain, while, if a ring R contains the rational numbers and has no non-zero prime D-ideals for a set D of derivations of R, then R is a D-simple ring (cf. [15], Corollary 1.5). Seidenberg proved that, if R is a domain which is a finitely generated algebra over a field, then R is a Der R-simple ring (where Der R denotes the set of all derivations of R), if, and only if, R is a regular ring (cf. [28], Theorems 3 and 5). Thus, if R is D-simple for any set D of derivations of R, then R is a regular ring. Hart [6] proved that this is actually true for the wider class of G-rings, which contains all finitely generated algebras and all complete local rings over fields and is closed under localization ([21], pp. 249-257). The converse of this statement is not true in general; e.g. the coordinate ring R of the real sphere (example 3.1 above), although it is regular, admits no derivation d such that R is a d-simple ring (cf [7], example iii). However, if R is a regular local ring of finitely generated type over a field k

(i.e. a localization of a finitely generated k-algebra at a prime ideal of R), Hart [7] constructed a derivation d of R, such that R is a d-simple ring. The above results show that for a wide class of rings of characteristic zero the d-simplicity is connected with the regularity, a property which requires from a ring R to have a "special" kind of dimension. In fact, in this case every maximal ideal M of R can be generated by dim $R_M$ elements, where $R_M$ denotes the localization of R at M. In particular, if R is a finitely generated algebra, then M can be generated by dim R elements.

Consider now the localization $R_P$ of a ring R at a proper prime ideal P of R. Then every derivation d of R extends to a derivation of $R_P$ by the usual rule of differentials for rational functions. We prove the following result:

**Theorem 1.6:** Let R and $R_P$ be as above. Then, if R is a D-simple ring with respect to a set D of derivations of R, $R_P$ is also D-simple. Conversely, if $R_P$ is D-simple with respect to a set D of derivations of $R_P$ whose restrictions to R are derivations of R, then R is also a D-simple ring.

*Proof:* Assume first that R is a D-simple ring, while $R_P$ is not. Let q be a prime ideal of $R_P$, such that $d(q) \subseteq q$ for all d in D. Then there exists a prime ideal Q of R not meeting P, such that $q = Q^e$

$$= \{ \frac{r}{s} : r \in Q, \ s \in R\text{-}P\} \text{ (cf. [2], Proposition 3.11). Since R}$$

is a D-simple ring, we can find r in Q, such that $d(r) \notin Q$

for all d in D. But $d(\dfrac{r}{1}) = d(r)$ is in

$Q^e = q$, therefore d(r) is in Q, which is absurd.

For the converse, assume that there exists a non zero prime ideal, say Q, of R, such that $d(Q) \subseteq Q$, for all derivations d of D. Then, for all r in Q, $d(r) \in Q \iff sd(r) \in Q$, for all s in P-R

$$\iff \frac{sd(r) - rd(s)}{s^2} = d(\frac{r}{s}) \in Q^e. \text{ Thus } R_P \text{ is not a D-}$$

simple ring, which is a contradiction. –

**Corollary 1.7:** Let R, $R_P$ and D be as in the converse statement of the above theorem. Assume further that R is a semi-local finitely generated algebra over a field k. Then R is a field.

*Proof:* By Noether's normalization Lemma R is an integral extension of a polynomial ring, say $k[x_1, x_2, \ldots, x_m]$, and dim R=m. Further, by the correspondence of maximal ideals in integral extensions, the number of maximal ideals of $k[x_1, x_2, \ldots, x_m]$ should be finite.

But $(x_1-a_1, x_2-a_2, \ldots, x_m-a_m)$ is a maximal ideal of $k[x_1, x_2, \ldots, x_m]$, for all choices of the $a_i$'s in k, i=1,2,…,m. Since k is of characteristic zero, it is an infinite field, therefore, if $m \neq 0$, then $k[x_1, x_2, \ldots, x_m]$ has an infinite number of maximal ideals, which is a contradiction. Thus m=0 and the Krull dimension of R is zero. But, by Theorem 1.6, R is a D-simple ring, therefore R is a domain and the result follows.-

The following example illustrates Theorem 1.6:

**Example 1.8:** Let R=$k[x_1,x_2]$ be a polynomial ring over a field k; then it is easy to check that $P=(x_1)$ is a prime ideal of R. By example 1.1 R is a D-simple ring, with D={ $\dfrac{\partial}{\partial x_1}$, $\dfrac{\partial}{\partial x_2}$ }. Therefore, by Theorem 1.6, $R_P$ is also D-simple, where D in this case denotes the set of extensions of $\dfrac{\partial}{\partial x_1}$, $\dfrac{\partial}{\partial x_2}$ to $R_P$.

## II.    FINITELY    GENERATED    D-SIMPLE ALGEBRAS OF DIMENSION 1

Let D be a non singleton set of derivations of a ring R, and assume that R is d-simple for some d in D, then obviously R is also a D-simple ring.

The converse is not true; e.g., although there exists a set D of 2 derivations of the coordinate ring R of the real sphere such that R is D-simple (cf. Example 1.3), there is no derivation d of R such that R is d-simple.(cf. [7], example iii).

For a second counter example notice that, if a complete local ring R is d-simple for some derivation d of R, then the Krull dimension of R is 1 (cf. [37], Theorem 2.3). Thus, although the power series ring R=k$[[x_1,x_2,…,x_n]]$ over a field k is { $\dfrac{\partial}{\partial x_1}$, $\dfrac{\partial}{\partial x_2}$,.... $\dfrac{\partial}{\partial x_n}$ }-simple (cf. Example 1.2), there is no derivation d of R such that R is a d-simple ring. Thus the "strongest" property of a ring R concerning its differential simplicity is to admit a derivation d, such that R is d-simple.

Notice that, if R is a differentially simple ring (and therefore a domain) of dimension 0, then R is a field. On the contrary, they are known several non trivial cases of d-simple rings of dimension 1, such as the polynomial and the power series rings in one variable over a field (see examples 1.1 and 1.2 above), as well as non trivial examples of 1-dimensional rings which are not d-simple for any derivation d (e.g. see example of [16]). Therefore it looks interesting to a search for necessary and sufficient conditions for the d-simplicity of 1-dimensional rings. We shall do this below for finitely generated algebras over a field k. In this case it is easy to prove:

**Proposition 2.1:** Let R be a finitely generated algebra over a field k, say R=$k[y_1,y_2,...,y_n]$. Assume further that there exists a derivation d of R such that R is a d-simple ring. Then R=$(d(y_1),d(y_2),...,d(y_n))$.
*Proof:* Theorem 2.4 of [37].-
In Theorem 2.4 of [37] was quoted without proof that the converse of Proposition 2.1 is also true. In order to prove this we need the following two lemmas:

**Lemma 2.2:** Let R and S be algebras over a field k such that S is integral over R, and let d be a derivation of R which extends to a derivation of S. Then, if P is a prime d-ideal of R and T is a prime ideal of S such that $T \cap R = P$, T is a d-deal of S.

*Proof:* Let $t$ be an element of T such that d($t$) is not in T. Since S is integral over R, there exists a positive integer n such that

$f(t)=t^n+a_{n-1}t^{n-1}+....+a_2t^2+a_1t+a_0=0$,

with $a_i$ in R for each i=0,1,...,n-1. Thus

$a_0=-t(t^{n-1}+a_{n-1}t^{n-2}+.....+a_2t+a_1)$

is in $T \cap R = P$, therefore d($a_0$) is also in P. But $df(t)=nt^{n-1}d(t)+d(a_{n-1})t^{n-1}+(n-1)a_{n-1}t^{n-2}d(t)+....$ $+d(a_2)t^2+2a_2td(t)+d(a_1)t+a_1d(t)+d(a_0)=0$.

Therefore we can write $a_1d(t)=At-d(a_0)$, with A in S. Thus $a_1d(t)$ is in T and therefore, since $d(t)$ is not in T, $a_1$ is in P. In the same way from $d^2[f(t)]=0$ we get that $a_2$ is in P and we keep going in the same way until we finally find that $a_i$ is in P for each i=0,1,...,n-1. Then $d^n[f(t)]=0$ implies that $n![d(t)]^n$ is in T, therefore $d(t)$ is in T, which contradicts our hypothesis.-

Using the previous lemma we can also prove:

**Lemma 2.3:** Let R be as in Proposition 3.1, and let d be a k-derivation of R. Then, if M is a maximal ideal of R invariant under d, d($y_i$) belongs to M for all i=1,2,...,n.
*Proof:* Consider the polynomial ring A=$k[x_1,x_2,...,x_n]$, then R $\cong$ A/I, where I is an ideal of A. Let f be natural map from A to R, then there exists a maximal ideal M' of A (containing I) such that M=f(M') and d induces a k-derivation of A (denoted also by d) such that d(M')$\subseteq$ M'.

In fact, if d($y_i$)=$h_i$, choose fixed $h_i'$ in A such that $f(h_i')=h_i$ and define d by $d(x_i)=h_i'$, i=1,...,n. Then, if m' is in M', we have that

$$d(m')=\sum_{i=1}^{n}\frac{\partial m'}{\partial x_i}h_i' , \text{ while } d[f(m')]=\sum_{i=1}^{n}\frac{\partial f(m')}{\partial y_i}h_i .$$

Thus $f[d(m')]=d[f(m')]$. But d(M)$\subseteq$ M, therefore $d[f(m')]$ is in M, i.e. $f[d(m')]$ is in M, which shows that $d(m')$ is in M'.
Next denote by K the integral closure of k and consider an arbitrary element, say F, of the polynomial ring B=$K[x_1,x_2,...,x_n]$. Let $a_1,...,a_j$ be the non zero coefficients of F, then, since $K'=k(a_1,...,a_j)$ is a finite extension of k, the polynomial ring $K'[x_1,x_2,...,x_n]$ (where F belongs) is a finitely generated A-module and therefore it is integral over A (cf. [13], Theorem17). Thus B is integral over A, therefore there exists a maximal ideal, say N, of B such that N $\cap$ A=M'(cf. [13], Theorem44 and [2], Corollary5.8, p.61). But d extends to a K-derivation of B, therefore by Lemma 2.1 we get that $d(N)\subseteq N$. Furthermore, by the Hilbert's Nullstellensatz

(zero-point Theorem), $N=(x_1-c_1,...,x_n-c_n)$ with $c_i$ in K for each i=1,2,...,n. Thus $d(x_i)$ is in N and therefore $d(y_i)$ is in M for all i=1,2,...,n.
We are ready now to prove:

**Theorem 2.4:** Let R and d be as in Lemma 3.3. Assume further that the Krull dimension of R is 1.Then R is a d-simple ring, if, and only if,
R =$(d(y_1),d(y_2),...,d(y_n))$.
*Proof:* Assume that R=$(d(y_1),d(y_2),...,d(y_n))$ and let M be a prime (and therefore maximal) ideal of R such that $d(M) \subseteq M$. Then, by Lemma 2.3, $d(y_i)$ belongs to M for all i=1,2,…,n, fact which contradicts our hypothesis. The converse is also true by Proposition 2.1. -
The following example illustrates the above theorem:

**Example 2.5:** The coordinate ring R= $\dfrac{k[x_1, x_2]}{(x_1^2 + x_2^2 - 1)}$

of the circle defined over a field k admits k-derivations d, such that R is a k-simple ring.
*Proof:* Theorem 2.2 of [39].-

**Remarks 2.6: (i)** As it becomes evident from the last part of the proof of Lemma 2.2, if k is an algebraically closed field, then Theorem 2.3 is a straightforward consequence of the Hilbert's Nullstellensatz. Also, if we take k to be the field of rational numbers, then, since every derivation of R is a k-derivation, Theorem 2.3 holds for any derivation of R.

**(ii)** A derivation d of a ring R is called a *simple derivation*, if R is a d-simple ring. They are certainly known examples of finitely generated algebras of dimension greater than 1, and even of infinite dimension, admitting simple derivations, typified by polynomial rings in finitely and infinitely many variables over a field k, and by Laurent polynomial rings in finitely many variables, say n, over k provided that the dimension of k (as a vector space) over the field of the rational numbers is greater or equal to n (cf. [37], section 3). In particular, for the polynomial ring $k[x,y]$ most of the published examples of its simple derivations with $d(x)=1$ are of the form d= $\dfrac{\partial}{\partial x} + F(x, y) \dfrac{\partial}{\partial y}$, where $F(x,y)$ is a polynomial of $k[x,y]$ with $deg_y F(x,y) \le 2$ (e.g. [3], [12], [18], etc) .
Nowicki [23] proved recently that $\dfrac{\partial}{\partial x} + (y^s + px) \dfrac{\partial}{\partial y}$,

where $s$ is an arbitrary positive integer and $0 \ne p$ in k, is a simple derivation of $k[x,y]$. His proof is based on the well known fact that a derivation d of $k[x,y]$ such that $d(x)=1$ is simple, if, and only if, d has no Darboux polynomials (cf. [22]).

**(iii)** It is well known that, if the coordinate ring of a variety , say Y, over a field k is regular, then Y is a smooth variety (i.e. Y has no singular points). This result, combined to the fact that a d-simple finitely generated algebra is a regular ring (cf. [28]), shows that, if R is the coordinate ring of a singular variety, then R admits no simple derivations. In [39] we have presented some characteristic examples of smooth varieties over k (e.g. cylinder, real torus considered as a 2-dimensional surface in 4 dimensions, etc), which admit at least one simple derivation. We emphasize that this is not true in general, the typical counter example being the coordinate ring of the real sphere (see example 1.3 above).

**(iv)** For some interesting properties of the differential ideals of a ring the reader may look [34].

## III.    ON A SPECIAL CLASS OF ITERATED SKEW POLYNOMIAL RINGS OF DERIVATION TYPE.

Let R be a ring (not necessarily commutative) with identity. In this section we are going to construct a special class of iterated skew polynomial rings defined with respect to any finite set D of derivations of R commuting to each other. The important thing with the skew polynomial rings of this class is that, if R is a commutative ring, then its D-simplicity is the necessary and sufficient condition for the simplicity of the corresponding skew polynomial ring.
In order to define skew polynomial rings over R, we extend the notion of a derivation of R as follows:

**Definition 3.1:** Let $f$ and $g$ be any two endomorphisms of R. Then a map $d$: R →R is called a *(f,g)-derivation* of R if $d(a+b)=d(a)+d(b)$ and
$d(ab)= g(a)d(b)+d(a)f(b)$, for all a,b in R.-
We have studied such kind of derivations of a ring R in [36]. Notice that, if both $f$ and $g$ coincide with the identity automorphism of R, then $d$ becomes a derivation of R. Further, if only $g$ is the identity automorphism of R, then $d$ is called a *f-derivation*, or *skew derivation* of R defined with respect to $f$.
Let now $f$ be a monomorphism of R and let $d$ be a f-derivation of R. Then we define a skew polynomial ring over R as follows:

**Definition 3.2:** Consider the set S of all polynomials in one variable, say x, over R. Define addition in S in the usual way and multiplication by the rule $xr=f(r)x+d(r)$ (1) and the distributive law of multiplication with respect to addition. Then S becomes a ring, called a *skew* (or *twisted) polynomial ring* over R and denoted by R[$x$; *f, d*].-
Such rings have been firstly introduced by Ore [24] in 1933 over a division ring R and that is why they are also known as *Ore extensions*. Initially they were used as counter examples, but eventually a broad bibliography was developed about them, since it was realized that they appear an important theoretical interest.
Notice that in the above definition we asked for f to be a monomorphism, and not a simple endomorphism of R, just to block the case for x to be a zero divisor of S. In fact, if there were $r_1 \ne r_2$ in R such that $f(r_1)=f(r_2)$, then
$x(r_1-r_2)=f(r_1)x+d(r_1)-[f(r_2)x+d(r_2)]=d(r_1-r_2)$,
therefore we should have $x(r_1-r_2)=0$, if $r_1-r_2$ is in the kernel of d.
If f is the identity automorphism of R, then S is denoted by R[$x$; *d*] and is called a skew polynomial ring *of derivation type* over R. Then relation (1) gives that xr=rx+d(r), for all r in R. Using this equation and applying induction on n one finds that

$$x^n r = \sum_{i=0}^{n} \binom{n}{i} d^i(r) x^{n-i},$$ for all r in R and all positive

integers n.
 The following example illustrates this case:

**Example 3.3:** Let $T[x_1]$ be a polynomial ring over a ring $T$, then the skew polynomial ring

$T[x_1][x_2; \dfrac{\partial}{\partial x_1}]$ of derivation type over $T[x_1]$ is called the

*first Weyl algebra over T* and is denoted by $A_1(T)$. It becomes evident that the elements of $A_1(T)$ are polynomials in two variables $x_1$ and $x_2$ over $T$, while

multiplication is defined by $x_1 t = t x_1$, $x_2 t = t x_2 + \dfrac{\partial t}{\partial x_1} = t x_2$

for all $t$ in $T$, $x_2 x_1 = x_1 x_2 + \dfrac{\partial x_1}{\partial x_1} =$

$= x_1 x_2 + 1$ and the distributive law.-

Notice also that, if d is the zero derivation of R, then the skew polynomial ring S is denoted by R[x; f], where multiplication is defined by $xr = f(r)x$ and the distributive law. Using the previous equation and applying induction on n one finds that $x^n r = f^n(r)x^n$, for all r in R and all positive integers n.

Next, following [5; p. xx, 14], we define skew polynomial rings in finitely many variables over R as follows:

**Definition 3.4:** Let $S_1 = R[x_1; f_1, d_1]$ be a skew polynomial ring over the ring $R$, where $f_1$ is a monomorphism and $d_1$ is a $f_1$-derivation of $R$. Then, if $f_2$ is a monomorphism and $d_2$ is a $f_2$-derivation of $S_1$, the skew polynomial ring $S_2 = S_1[x_2; f_2, d_2]$ is called an *iterated skew polynomial ring* over R and is denoted by $S_2 = R[x_1; f_1, d_1][x_2; f_2, d_2]$. Keeping the same notation let us consider the finite sets $H = \{f_1, f_2, ..., f_n\}$, where $f_i$ is a monomorphism of $S_{i-1}$ and $D = \{d_1, d_2, ... d_n\}$, where $d_i$ is a $f_i$-derivation of $S_{i-1}$, i=1,0,...,n (we set $S_0 = R$). Then by induction we define the iterated skew polynomial ring $S_n = R[x_1; f_1, d_1][x_2; f_2, d_2] ..... [x_n; f_n, d_n]$ (in n indeterminates) over R . In order to simplify our notation we shall denote this ring by $S_n = R[x; H, D]$.

If all the elements of $H$ are identity automorphisms, then $S_n$ is denoted by R[x; D] (iterated skew polynomial ring *of derivation type* over R), while if all the elements of D are zero derivations, then we have the iterated skew polynomial ring $S_n = R[x; H]$.

**Examples 3.5: (i)** The first Weyl algebra $A_1(T)$ over a ring $T$ (cf. example 3.3) is an iterated skew polynomial ring of derivation type in two variables over $T$ of the form

$T[x_1; d][x_2, \dfrac{\partial}{\partial x_1}]$, where $d$ denotes the zero derivation of

T (in this case

$T[x_1; d]$ is the ordinary polynomial ring $T[x_1]$).

**(ii)** Set $R = A_1(T)$. Then the first Weyl algebra $A_1(R)$ over $R$ is called the *second Weyl algebra over T* and is denoted by $A_2(T)$. Obviously we have that

$A_2(T) = A_1[A_1(T)] = T[x_1][x_2; \dfrac{\partial}{\partial x_1}][x_3; \dfrac{\partial}{\partial x_2}]$.

**(iii)** Consider the set of all polynomials in n+1 variables, say $x_1, x_2, ..., x_n, x_{n+1}$, over a ring $T$. Then the *n-th Weyl algebra $A_n(T)$ over $T$* is defined by induction as $A_n(T) = A_1[A_{n-1}(T)]$. Obviously we have that

$A_n(T) = T[x_1][x_2; \dfrac{\partial}{\partial x_1}][x_3; \dfrac{\partial}{\partial x_2}] ...... [x_{n+1}; \dfrac{\partial}{\partial x_n}] =$

$= T[x; D]$, with $D = \{d, \dfrac{\partial}{\partial x_1}, \dfrac{\partial}{\partial x_2} ......, \dfrac{\partial}{\partial x_n}\}$, where $d$

denotes the zero derivation of T.-

Next, given a finite set $D$ of derivations of R <u>commuting</u> to each other, we shall construct an iterated skew polynomial ring of derivation type over R of the form $R[x; D]$. For this, we need the following lemma:

**Lemma 3.6:** Let $R$ be a ring, let $d$ be a derivation of $R$, let $S = R[x; d]$ be the skew polynomial ring of derivation type over $R$ defined with respect to $d$, and let $d_1$ be another derivation of $R$. Then $d_1$ can be extended to a derivation of $S$ by $d_1(x) = 0$, if and only if $d_1$ commutes with $d$.
*Proof:* Assume first that $d_1$ is extended to a derivation of S by $d_1(x) = 0$. Then, given r in R, we have that $xd_1(r) = d_1(xr) = d_1[rx + d(r)] = d_1(rx) + d_1[d(r)] = d_1(r)x + d_1[d(r)]$ (2).
But $xd_1(r) = d_1(r)x + d[d_1(r)]$ (3) and therefore (2) implies that $d_1[d(r)] = d[d_1(r)]$ for all r in R, or $d_1 \circ d = = d \circ d_1$.
Conversely assume that $d_1$ commutes with $d$. Then $d_1$ could be extended to a derivation of S, if $d_1(x)$ could be defined in a way compatible with multiplication in S. But, given r in R, we have by (2)
that $d_1(xr) = d_1(r)x + d_1[d(r)]$. Therefore, since $d_1 \circ d = = d \circ d_1$, (3) gives that $d_1(xr) = xd_1(r)$ (4).
But $d_1(xr) = d_1(x)r + xd_1(r)$, therefore we must have $d_1(x)r = 0$ for all r in R, which is true if we extend $d_1$ to S by $d_1(x) = 0$.-
We are ready now to prove the following theorem:

**Theorem 3.7:** Let $R$ be a ring, and let $D = \{d_1, d_2, ...., d_n\}$ be a finite set of derivations of $R$.
Denote by $S_i$ the set of all polynomials in indeterminates, say $x_1, x_2, ..., x_i$ with coefficients in $R$, for i=0,1,2,....,n, where $S_0 = R$. Define in $S_n$ addition in the usual way and multiplication by the rules $x_i r = r x_i + d_i(r)$ for all $r$ in R, $x_i x_j = x_j x_i$ (5) for i,j=1,2,....n and the distributive law of multiplication with respect to addition.
Then $S_i$ is a skew polynomial ring (of derivation type) over $S_{i-1}$, for all i=1,2,...n, if and only if $d_i \circ d_j = d_j \circ d_i$, for all i,j=1,2,....n.
*Proof:* Assume first that the elements of $D$ commute to each other. Obviously $S_1 = R[x_1; d_1]$ is a skew polynomial ring over R. We apply induction on n. In fact, assume that $S_i$ is a skew polynomial ring over $S_{i-1}$ for each $i \le n-1$. By Lemma 3.6 $d_n$ is extended to a derivation of $S_1$ by $d_n(x_1) = 0$. But, by our inductive hypothesis, $S_2 = S_1[x_2, d_2]$ is a skew polynomial ring over $S_1$. Therefore, by Lemma 3.6 again, $d_n$, being a derivation of $S_1$, is extended to a derivation of $S_2$ by $d_n(x_2) = 0$ and so on. We keep going in the same way until we find, after a finite number of steps, that $d_n$ is extended to a derivation of $S_{n-1}$ by $d_n(x_i) = 0$,

i=1,2,…,n-1.                    Let                    now

$h = \sum\limits_{a_1+a_2+...+a_{n-1}=0}^{k} r_{a_1,a_2,...a_{n-1}} x_1^{a_1} x_2^{a_2} ... x_{n-1}^{a_{n-1}}$ be an element

of $S_{n-1}$, with $r_{a_1,a_2,......,a_{n-1}}$ in $R$. For reasons of brevity we write $h = \sum\limits_{(a)} r^{(a)} x^{(a)}$ .

Then

$x_n h = x_n (\sum\limits_{(a)} r^{(a)} x^{(a)}) = \sum\limits_{(a)} x_n r^{(a)} x^{(a)} =$

$\sum\limits_{(a)} [r^{(a)} x_n + d_n(r^{(a)})] x^{(a)} =$

$(\sum\limits_{(a)} r^{(a)} x^{(a)}) x_n + \sum\limits_{(a)} d_n(r^{(a)}) x^{(a)} =$

$= h x_n + d_n(h)$. Thus $S_n = S_{n-1}[x_n; d_n]$ is a skew polynomial ring over $S_{n-1}$.

Conversely, assume that $S_i$ is a skew polynomial ring over $S_{i-1}$, for all i=1,2,..,n. Then, since

$S_n = S_{n-1}[x_n, d_n]$ is a skew polynomial ring over $S_{n-1}$ and $x_i$ belongs to $S_{n-1}$ for each i, $1 \leq i < n$, we have that $x_n x_i = x_i x_n + d_n(x_i)$, therefore by the second of rules (5) we get that $d_n(x_i) = 0$.

Further, given r in R we have that $x_i r = r x_i + d_i(r)$, therefore $d_n(x_i r) = d_n(r x_i) + d_n[d_i(r)]$, or $x_i d_n(r) = d_n(r) x_i + d_n[d_i(r)]$.

But $x_i d_n(r) = d_n(r) x_i + d_i[d_n(r)]$ and the last two equalities imply that $d_n o d_i = d_i o d_n$. In the same way and for each j=1,2,…,n one can show that $d_j o d_i = d_i o d_j$, for all i, $1 \leq i < j$ and this completes the proof.-

The above theorem shows that, if the derivations of $D$ commute to each other, then $S_n = R[x; D]$, with addition defined in the usual way and multiplication by the rules (5) and the distributive law, is an iterated skew polynomial ring of derivation type over R.

A more general version of this theorem, concerning the case of an iterated skew polynomial ring over R of the form $R[x; H, D]$ with $H$ a finite set of monomorphisms of $R$ and $D$ a finite set of skew derivations of $R$ defined with respect to the monomorphisms of $H$, was quoted without an explicit proof in an earlier paper of the author (see Theorem 2.4 of [30]).

The following example illustrates Theorem 3.7:

**Example 3.8:** Let $R = T[y_1, y_2, ......, y_n]$ be a polynomial ring over a ring $T$. Set $D = \{\frac{\partial}{\partial y_1}, \frac{\partial}{\partial y_2}, ...., \frac{\partial}{\partial y_n}\}$. Since the elements of $R$ are polynomials with coefficients in $T$, their partial derivatives are continuous functions. Therefore, by the Young's classical theorem of differential calculus, the derivations of $D$ commute to each other.

Consider the set $S_n$ of all polynomials in n indeterminates, say $x_1, x_2, ....., x_n$, with coefficients in $R$ and define in it addition in the usual way and multiplication by the rules

$h x_i = x_i h + \frac{\partial h}{\partial y_i}$ for all $h$ in R and $x_i x_j = x_j x_i$, i,j=1,…,n. Then,

by Theorem 3.7, $S_n = R[x; D]$ is an iterated skew polynomial ring of derivation type over the polynomial ring $R$. It is easy to check that $y_i x_i = x_i y_i + 1$, while $y_i x_j = x_j y_i$ for i $\neq$ j ,

i, j=1,2,…,n.-

Notice that, in order to have $x_i x_j = x_j x_i$, it is necessary to have $d_i o d_j = d_j o d_i$, for all i,j=1,2,…,n. In fact, given r in R, we have

$x_i x_j r = x_i [r x_j + d_j(r)] = (x_i r) x_j + x_i d_j(r) =$

$= r x_i x_j + d_i(r) x_j + d_j(r) x_i + (d_i o d_j)(r)$   (6).

In the same way one finds that

$x_j x_i r = r x_j x_i + d_j(r) x_i + d_i(r) x_j + (d_j o d_i)(r)$   (7)

and the result follows by equating the right hand sides of (6) and (7).

However, if there exist x's in $S_n$ not commuting to each other, this does not mean that the same must happen with the derivations of $D$. An example illustrating this situation is the n-th Weyl algebra $A_n(T)$ over a ring $T$, $n \geq 1$ (see examples 3.5). In fact, the elements of $A_n(T)$ are polynomials over $T$ and therefore, as in example 3.8, the derivations of $D$ commute to each other, while the x's they don't commute (see example 3.3).

Skew polynomial rings (of derivation type) in finitely many indeterminates (not necessarily commuting to each other) over a ring R with respect to a finite set of derivations of R have been firstly introduced by Kishimoto [14] . In these rings, which are not necessarily iterated skew polynomial rings over R, multiplication is defined by the first only of rules (5) and the distributive law.

We proceed by recalling the following definition:

**3.9 Definition:** A derivation $d$ of a ring $R$ is called an *inner derivation* induced by an element $s$ of R, if

$d(r) = sr - rs$, for all $r$ in R. A derivation of $R$ which is not inner is called an *outer derivation*.-

Obviously, if $R$ is a commutative ring and $d$ is an inner derivation of $R$, then $d = 0$.

Notice also that an outer derivation of $R$ is possible to be extendable to an inner derivation of a skew polynomial ring of derivation type over R. For example, $d = \frac{\partial}{\partial x_1}$ is an outer derivation of the polynomial ring R=$T[x_1]$ over a ring $T$, because $x_1$ is a central element of R and $d(x_1) = 1 \neq 0$. But $d$ is extended to an inner derivation of $A_1(T) = R[x_2, d]$ induced by $x_2$, because, by the definition of multiplication in $A_1(T)$, we have that $x_2 f = f x_2 + d(f)$, or $d(f) = x_2 f - f x_2$, for all $f$ in R.

For the general case we prove the following result:

**3.10 Proposition:** Let $d$ and $d'$ be two derivations of a ring $R$. Assume that $d'$ is extended to an inner derivation of the skew polynomial ring $S = R[x; d]$ induced by an

element $f = \sum\limits_{i=0}^{n} a_i x^i$ of S. Then $a_n$ is a central element of $R$

and

$ra_k - a_k r = \sum_{i=k+1}^{n} a_i \binom{i}{i-k} d^{i-k}$ (r), for all $r$ in R and each k=0,1,….,n-1.

*Proof:* For all $r$ in R, $d'(r) = fr - rf =$

$\sum_{i=0}^{n} [a_i(x^i\, r) - ra_i x^i]$ is in R, where

$x^i r = \sum_{j=0}^{i} \binom{i}{j} d^j(r) x^{i-j}$ and the result follows by straightforward calculations.-

However, if $f$ is a monomorphism of $R$ and $d$ is an outer derivation of $R$ extendable to a derivation of the skew polynomial ring $S=R[x;\, f]$, then $d$ cannot be an inner derivation of $S$. In fact, assume that $d$ is an inner derivation of $S$ induced by an element

$h = \sum_{i=0}^{n} a_i x^i$ of $S$. Then, for all $r$ in R $d(r) = hr - rh$

$= \sum_{i=0}^{n} [a_i(x^i\, r) - ra_i x^i] = \sum_{i=0}^{n} [a_i f^i(r) - ra_i] x^i$ must be in R,

wherefrom we get that $d(r) = a_0 r - ra_0$ for all $r$ in R. Thus $d$ is an inner derivation $R$ induced by $a_0$, which contradicts our hypothesis.

We are ready now to state the following important remarks:

**3.11 Remarks: (i)** Let $R$ be a ring of characteristic zero and let $D = \{d_1, …….., d_n\}$ be a finite set of derivations of R commuting to each other. Assume further that $R$ is a $D$-simple ring and that $d_i$ is an outer derivation of $S_{i-1}$, i=1,2,…,n (where $S_0=R$). Then $S_n=R[x;\, D]$ is a simple ring [29; Theorem 3.4].

Conversely, if $S_n=R[x;\, D]$ is a simple ring, then no element of $D$ is an inner derivation of $R$ induced by an element of $\bigcap_{d \in D} Kerd$ (where $Kerd$ denotes the kernel of $d$) and $R$ is a $D$-simple ring [29; Theorem 3.3].

**(ii)** Analogous results hold (although their statement is a little bit more complicated), if $R$ is of prime characteristic [35; Theorems 2.3 and 2.4].

**(iii)** As an immediate consequence of the previous remarks, if $R$ is a commutative ring and $D$ is a finite set of derivations of $R$ commuting to each other, then $S_n=R[x;\, D]$ is a simple ring, if and only if $R$ is a $D$-simple ring.

**(iv)** Let $S_n^*=R[x;\, D]$ be a skew polynomial ring in finitely many indeterminates (not necessarily commuting to each other under multiplication) over a ring $R$ of characteristic zero, defined with respect to a set $D=\{d_1, d_2, …, d_n\}$ of derivations of $R$ (we recall that multiplication in $S_n^*$ is defined by the first of rules (5) only and the distributive law). Then, if $R$ is a $D$-simple ring and $d_i(C(S_{i-1}^*) \cap R) \neq 0$, where $C(S_{i-1}^*)$ denotes the center of $S_{i-1}^*$ and $S_0^*=R$, for all i=1,2,…,n , $S_n^*$ is a simple ring [38; Theorem 2.1].

Obviously, if $d_i(C(S_{i-1}^*) \cap R) \neq 0$, then $d_i$ is an outer derivation of $S_{i-1}^*$. Therefore, if the elements of $D$ commute to each other, the above result is a weaker form of the result mentioned in remark (i). Also, if $R$ is a commutative ring, then the

$D$-simplicity of R implies the simplicity of $S_n^*$.

**(v)** Let $D$ be a set of derivations of a ring $R$. Then, we call a $D$- ideal $I$ of $R$ a $D$-prime ideal, if, given any two $D$-ideals $A$, $B$ of $R$ such that $AB \subseteq I$, we have either $A \subseteq I$, or $B \subseteq I$.

There is a connection among the prime ideals of the iterated skew polynomial ring $S=R[x;D]$ of Theorem 3.7 and the $D$-prime ideals of $R$. In fact, if $P$ is a prime ideal of $S$, then $P \cap S$ is a $D$-prime ideal of $R$ [32; Theorem 2.2]. Conversely, if $I$ is a $D$-prime ideal of $R$, then $IS$ is a prime ideal of S [32; Theorem 2.6].

**(vi)** Analogous relations hold among the semiprime ideals of $S$ and the $D$-semiprime ideals of $R$ [33]. We recall that a $D$-ideal $I$ of $R$ is called $D$-semiprime ideal, if for all $D$-ideals $A$ of $R$ such that $A^k \subseteq I$, with $k$ a positive integer, we have that $A \subseteq I$.

**(vii)** Let $H$ be a finite set of automorphisms of a ring $R$ commuting to each other. Then the iterated skew polynomial ring $S=R[x;\, H]$ of Theorem 3.7 is not simple, since $x_1^{a_2} x_2^{a_2} ….. x_n^{a_n} S$, with $a_1$, $a_2$, ..., $a_n$ non negative integers not all zero, is obviously a non zero, proper ideal of $R$.

The quotient ring $T$ of $S$ with respect to the Ore subset $C=\{\ x_1^{a_2} x_2^{a_2} ….. x_n^{a_n}:\ a_1,\ a_2,\ …,\ a_n \in \mathbf{N}\}$, where $\mathbf{N}$ is the set of non negative integers, is denoted by $T=R[x, x^{-1};\, H]$ and is called an *iterated skew Laurent polynomial ring* (in n indeterminates) over $R$. In [31] we have obtained the necessary and sufficient conditions for the simplicity of $T$.

## IV.  CONCLUSIONS AND OPEN PROBLEMS

The following conclusions, together with some open questions for future research, can be drawn from the discussion presented in the paper:

- A differentially simple ring of prime characteristic has zero Krull dimension. This means that for non trivial cases it cannot be an integral domain (otherwise it is a field), while just the opposite happens for rings of characteristic zero.
- For rings of characteristic zero and apart from special cases (e.g. 1-dimensional finitely generated algebras over a field, regular local

rings of finitely generated type, etc) there is not known any general criterion to decide whether, or not, a given ring is differentially simple.

- Although it seems that in most cases there exists a connection between the differential simplicity and the Krull dimension of a ring, it looks difficult to find a criterion (for the case of characteristic zero) to be found in future, since there exist known examples of differentially simple rings even of infinite dimension, e.g. polynomial rings in infinitely many variables. On the contrary, it is possible to be found in future new examples (isolated, or categories) of differentially simple rings.

- Some times the differential simplicity of a ring can be connected with its geometric characteristics, e.g. coordinate rings of certain smooth varieties over a field of characteristic zero, like the circle, the real sphere, the cylinder, the real torus, etc. It is certainly known however that not all the coordinate rings of such varieties admit simple derivations (e.g. real sphere). But what happens for the smooth varieties over an algebraically closed field? Very possibly, for example, the coordinate ring

$$\frac{C[x,y,z]}{(x^2+y^2+z^2-1)}$$ of the complex sphere,

where $C$ denotes the field of complex numbers, admits a simple derivation (we have studied the unpublished work of a colleague, who claims that he has constructed such a derivation, although we have some doubts for the correctness of the proof given). This possibility raises the following question: *If R is the coordinate ring of a smooth variety over an algebraically closed field, does R admit simple derivations?* The answer to this question is not known (at least to us).

- In section 3 we have constructed a special class of iterated skew polynomial rings of derivation type over a ring R, such that, if R is commutative, then its differential simplicity is the necessary and sufficient condition for the simplicity of the corresponding skew polynomial ring S of the above class. There is also a connection among the prime ideals of S and the, so called by us, D-prime ideals of R.

- Not all known properties of skew polynomial rings of derivation type in one indeterminate over R can be transferred by induction to the iterated skew polynomial rings of the above class. For example, this happens with the result about the simplicity of S, which for n=1 is due to Jordan [11]. Therefore a future research activity could concern the effort to transfer other such properties from $S_1$ to $S_n$. An other interesting activity could also concern the effort to transfer these properties to skew polynomial rings in finitely many indeterminates over R, which are not commuting under multiplication (e.g. see remarks 3.11, iv).

## IV. EPILOGUE

It is well known that progress in mathematics is usually achieved in two ways.

The first of them starts from the effort of applying mathematics to solve real world problems, or practical problems of our everyday life. Sometimes however, the already known mathematical theories are proved to be not suitable (or not enough) for solving the corresponding problem. In such cases the researchers are forced in trying to invent new mathematical theories in order to achieve their purpose. Characteristic examples are the recently developed theories of Fuzzy Sets [40] and of Chaotic Dynamics and Fractals ([20], [26]), which they have found lots of applications in practice (e.g. [1], [8], [9], [25], etc).

The second way concerns the pure mathematical research, i.e. the process of "doing mathematics for mathematics''. But, the amazing thing is that frequently, theories developed through pure mathematical research, find, some times many years after their development, successful applications to real situations. The best example for this is probably the use by Einstein of the Rieman's Geometry as a tool for the development of his Relativity Theory.

This amazing phenomenon was fist underlined by the great Plato, who believed in an already existing "Universe of Mathematics", wherefrom people eventually find some elements. Many centuries later H. O. Pollak [27] presenting the "Circle of Modelling" between the real world and the universe of mathematics considered the *"Classical Applied Mathematics"* and the *"Applicable Mathematics"* as two intersected, but not equal sets, the latter being topics from mathematics with a great mathematical interest, but without any visible, for the moment, applications, although it is possible to find such applications in future (for more details and a characteristic diagram see [36 ; section 2]).

In our case, derivations and skew polynomial rings are topics of the modern abstract algebra, which is certainly a field of pure mathematical research. However skew polynomial rings have found recently two very interesting applications in real practice and that is why the researchers' interest has been renewed about them.

The former of these applications concerns the use of *quantum groups*, i.e. Hopf algebras having in addition a structure analogous to that of a Lee group [19], as a tool in Theoretical Physics. In fact, it has been observed that many of these algebras can be expressed in the form of an iterated skew polynomial ring.

The latter application concerns the use of skew polynomial rings in *coding theory*. More explicitly, Jatengaokar [10] studied in detail the structure of monomorphism skewed polynomial rings over semisimple rings and showed that they are indeed direct sums of matrix rings. The entries in these rings came from subrings of matrix rings over skew polynomial rings. Jatengaokar's structure theorem has been recently used in coding theory to analyze the structure of certain convoliutional codes and, for the same purpose, Jatengaokar's results have been extended to the case of general skew polynomial rings, i.e. polynomial rings that are skewed by both a monomorphism and a derivation [17].

**Michael Gr. Voskoglou** (B.Sc., M.Sc., M.Phil., Ph.D) is currently a Professor in School of Technological Applications of Graduate Technological Educational Institute (T. E. I.) of Patras, in Greece, where he teaches Mathematics and Operations' Research. He is author of 5 books on Mathematics and Operation's Research (in Greek language) and of about 200 papers published in scientific periodicals and proceedings of international conferences in U.S.A., Great Britain, Italy, Germany, Portugal, Serbia, Montenegro, Bulgaria, Poland, Hungary, Czech Republic, Ukraine, Turkey, China, India, Singapore, Cyprus and Greece,, with many references from other researchers. He is also reviewer of the AMS Mathematical Reviews and referee in several mathematical magazines and international conferences.

His research interests include Algebra (Ring Theory - Linear Algebra), Applications of Markov Chains and Fuzzy Sets to Management, Economics, Education and Artificial Intelligence and Didactics of Mathematics.